\documentclass[12pt]{article}
\usepackage{amsmath,amssymb,graphicx,amsthm}
\usepackage{subfigure}

\newtheorem{theorem}{Theorem}[section]
\newtheorem{lemma}[theorem]{Lemma}
\newtheorem{corollary}[theorem]{Corollary}
\newtheorem{problem}[theorem]{Problem}

\newtheorem{definition}[theorem]{Definition}

\newcommand\lref[1]{Lemma~\ref{lem:#1}}
\newcommand\tref[1]{Theorem~\ref{thm:#1}}
\newcommand\cref[1]{Corollary~\ref{cor:#1}}

\begin{document}

\title{A note about online nonrepetitive coloring $k$-trees}

\author{Bal\'azs Keszegh\thanks{Alfr\'ed R\'enyi Institute of Mathematics and MTA-ELTE Lend\"ulet Combinatorial Geometry Research Group. Research supported by the Lend\"ulet program of the Hungarian Academy of Sciences (MTA), under grant number LP2017-19/2017 and the National Research, Development and Innovation Office -- NKFIH under the grant K 116769.}\and Xuding Zhu\thanks{Department of Mathematics, Zhejiang Normal University. Research supported by NSFC 11971438.}}

\maketitle
\begin{abstract}
We prove that it is always possible to color online nonrepetitively any (partial) $k$-tree (that is, graphs with tree-width at most $k$) with $4^k$ colors. This implies that it is always possible to color online nonrepetitively cycles, trees and series-parallel graphs with $16$ colors. Our results generalize the respective (offline) nonrepetitive coloring results.
\end{abstract}

\section{Introduction}
A sequence $x_1...x_{2l}$ is a \textit{repetition} if $x_i=x_{k+i}$ holds for all $1 \le i \le l$. A sequence is \textit{nonrepetitive} if it does not contain a string of consecutive entries forming a repetition. In 1906, Thue \cite{T} found an infinite nonrepetitive sequence using only three symbols.

Alon, Grytzuk, Ha\l uszczak, Riordan \cite{AGHR} generalized the notion of nonrepetitiveness to graphs: a coloring $c$ of a graph $G$ is \textit{nonrepetitive } if there is no path $v_1,...,v_{2l}$ in $G$ such that the string $c(v_1),...,c(v_{2l})$ is a repetition. The Thue-coloring number of a graph $G$ is the least integer $\pi(G)$ such that there exists a nonrepetitive  coloring $c$ of $G$ with $c:V(G)\rightarrow \{1,2,...,\pi(G)\}$. With this notation Thue's result says $\pi(P_{\infty})=3$ (the fact that 2 colors are not enough can be easily seen even for a path of length 4). A survey and a nice introduction to the topic is \cite{G}.

In this paper we investigate an online variant of Thue's theorem, where in each step we get a new vertex $v$ together with some edges that connect $v$ to previous vertices. Deletion of any edge is also allowed in any step, as anyway this just makes easier to color. We have to color the new vertex such that the coloring of the new graph remains nonrepetitive. Note that we color a vertex immediately when it is coming and no recoloring of previous vertices is allowed.

We now define this coloring scheme formally. At the beginning we start with some graph $G_0$ (it is the empty graph if not defined otherwise) and in each step $t$, we get a new vertex $v_t$ and $D_t$, a set of edges connecting $v_t$ to some previous vertices $v_i, i<t$ and $C_t\subset E(G_{t-1})$, a set of edges in the previous graph, and we set $V(G_t)=V(G)\cup\{v_t\}$ and $E(G_t)=E(G)\cup D_t\setminus C_t$. In this step we have to color $v_t$ such that $G_t$ is nonrepetitively colored. If under some restriction $R$ on $D_t$ and $C_t$ we can continue this nonrepetitive coloring process until we reach the final graph $G_n$ ($n=\infty$ is allowed), irrespective of what exactly the $D_t$'s are, we say that we can {\em color nonrepetitively online} any graph that can be generated according to $R$.

We list some natural graph classes with the corresponding restriction on $D_t$ and $C_t$ ($G_0$ is always the graph with a single vertex if not stated otherwise):

\begin{definition}
{\em Online left-to-right path:}
 all time $G_{t-1}$ is a path and in a step we append a vertex at the same end of this path ($D_t=v_tv_{t-1}$ and $C_t=\emptyset$).

{\em Online path:}
all time $G_{t-1}$ is a path and in a step we either append a vertex at one end of this path or we subdivide an edge by adding the new vertex connecting two consecutive vertices on the path and deleting the edge between these two points.

{\em Online tree:} 
all time $G_{t-1}$ is a tree and in a step we either connect the new vertex to one of the vertices, or we subdivide an edge by adding the new vertex connecting two consecutive vertices of the tree and deleting the edge between these two points.

{\em Online cycle:}
$G_0$ is a triangle, all time $G_{t-1}$ is a cycle and in a step we subdivide an edge by adding the new vertex connecting two consecutive vertices on the cycle and deleting the edge between these two points.

{\em Online series-parallel graph:}
all time $G_{t-1}$ is a series-parallel graph and in a step we either connect the new vertex to one of the vertices, or we subdivide an edge by adding the new vertex connecting two adjacent vertices of the graph and deleting the edge between these two points, or we just add the new vertex connecting two adjacent vertices of the graph.

{\em Online partial k-tree:}
all time $G_{t-1}$ is a partial $k$-tree and in a step we connect the new vertex to a clique of size at most $k$ and delete some edges from $G_{t-1}$.

{\em Online k-tree:}
$G_0$ is a complete graph on $k+1$ vertices, all time $G_{t-1}$ is a $k$-tree and in a step we connect the new vertex to a clique of size $k$.
\end{definition}

Note that the left-to-right online version is uninteresting as an infinite (offline) Thue-coloring gives a coloring with $3$ colors in this setting too. However, the list-coloring version is non-trivial already of this most simple online graph, as we will see in the last section, where we extend our scope to online nonrepetitive list-colorings.

Online paths are a natural class, they can be imagined also as points on a line and in each step a new point is placed somewhere on the line. We need that at all times all intervals are non-repetitively colored. Many other types of colorings are investigated for this online set of points, see, eg., the intriguing problem about online conflict-free coloring points on the line \cite{fiatetal}.

Finishing the introduction, to put our online model in context, we briefly mention other online models that gained attention in recent years related to various colorings. Our \emph{online} model, where vertices arrive one-by-one is probably regarded most widely, there are many results for various colorings and for hypergraphs as well. In the \emph{semi-online} model a vertex does not need to be colored immediately. In the \emph{dynamic} (also called \emph{quasi-online}) model we know in advance the whole sequence of vertices and edges that will arrive while we still require that in each step the coloring has some required property. Finally, the \emph{online choosability} (also called \emph{paintability}) model is an online variant of list-coloring, where the whole graph is given in advance and the lists of the colors are arriving one-by-one: in each step for the next color a subset is given which can receive this color and one has to color an independent subset of this subset with this color.

\subsection{Our results}

\begin{theorem}\label{thm:onlinektree}
Using $4^k$ colors we can color online nonrepetitively any $k$-tree.
\end{theorem}

Coloring nonrepetitively a subgraph $G'$ of a graph $G$ always requires a subset of the conditions compared to nonrepetitively coloring the original graph $G$. In particular, for partial $k$-trees we can always suppose that $C_t=\emptyset$, i.e., we do not delete edges, and also we can suppose that $|D_t|=\min(k,|V_{t-1}|)$, i.e., we always connect $v_t$ to exactly $k$ previous vertices or to all previous vertices if there are less vertices. Thus to color online nonrepetitively partial $k$-trees it is enough to color online nonrepetitively (non-partial) $k$-trees. We note that (offline) partial $k$-trees are exactly the graphs that have tree-width bounded by $k$.

\begin{corollary}\label{thm:onlinepartialktree}
Using $4^k$ colors we can color online nonrepetitively any partial $k$-tree.
\end{corollary}

This strengthens the result of K\"undgen and Pelsmajer \cite{KP}, who showed that graphs that have tree-width bounded by $k$, that is, partial $k$-trees, are (offline) nonrepetitively $4^k$-colorable. The proof of our result is based on this offline result.

It is easy to see that online trees, cycles and series-parallel graphs are online partial $2$-trees and thus it follows that:

\begin{corollary}
Using $16$ colors we can color online nonrepetitively any tree, any cycle and any series-parallel graph.
\end{corollary}

For paths the upper bound can be slightly improved by examining the proof more carefully.

\begin{theorem}\label{thm:onlinepath}
Using $12$ colors we can color online nonrepetitively any path.
\end{theorem}

We note that after solving Theorem \ref{thm:onlinepath} it came to our knowledge that it  was proved simultaneously and independently by Grytczuk et al. \cite{thuev2} using the same methods. They also proved by a simple argument that $5$ is a lower bound. We include the short proof of Theorem \ref{thm:onlinepath} for the sake of completeness and because it raises some open problems which we want to highlight.

In the next section we prove these theorems and in the last section we propose the investigation of the list-coloring version of these problems and prove some preliminary results.

\section{Online nonrepetitive colorings}\label{sec:onlineproper}

We first prove \tref{onlinektree} which guarantees a nonrepetitive online coloring of $k$-trees using $4^k$ colors, then we prove \tref{onlinepath} that improves the needed number of colors from $16$ to $12$ in the special case of online coloring paths.

\begin{proof}[Proof of \tref{onlinektree}]
During the proof we consider the (online built) graphs $G_t$ with their vertices numbered according to the online process as $v_1,v_2\dots$, that is $V(G_t)=\{v_1,v_2\dots v_t\}$.

We first define an (infinite) universal graph $U$, which is a $k$-tree, and for any $G_t$ that can be built online (i.e., there exists a sequence of graphs $G_0\subset G_1\subset\dots G_t$ that follows the restriction of being an online $k$-tree), we give an injection $inj_{G_t}$ from the underlying numbered vertices $v_1,v_2,\dots v_t$ to the vertices of $U$ such that the image of $G_t$ is a subgraph of $U$. These injections will be incremental, i.e. if $G_t$ is an extension of $G_{t-1}$ following the rules of building an online $k$-tree, then $inj_{G_t}$ is an extension of $inj_{G_{t-1}}$, i.e. these two injections coincide on the first $t-1$ vertices.
Thus, if for two graphs $G$ and $G'$ the induced subgraphs on their first $t_0$ vertices are (numbered) isomorphic, then these first $t_0$ vertices are embedded into $U$ by $inj_{G}$ and $inj_{G'}$ in the same way.
 
Assume we have shown such a universal graph $U$ that it is also a $k$-tree. As $U$ is infinite, by this we mean that it is a union of an infinite sequence of increasing $k$-trees $U_1\subset U_2\subset \dots =U$, where in $U$ the vertices of $U_i$ induce the graph $U_i$. Now K\"undgen and Pelsmajer \cite{KP} showed that graphs that have tree-width bounded by $k$, thus in particular $k$-trees, are nonrepetitively $4^k$-colorable. Thus every $U_i$ is non-repetitively $4^k$-colorable. Using K\"onig's Lemma, the classical result about infinite sequences, we can deduce that $U$ is nonrepetitively $4^k$-colorable as well. By that we mean a coloring for which every finite path in $U$ is nonrepetitively colored. Let us regard this coloring as a function $c$ on the vertices of $U$.

The coloring process is the following. We maintain that at any time $t$, every vertex $v_i$ $(i\le t)$ of $G_t$ is colored to $c'(v_i)=c(inj_{G_t}(v_i))$. Suppose that we colored $G_1,\dots G_{t-1}$ and then we get the new vertex and the new edges incident to it, together with $G_{t-1}$ they form $G_t$. Now by definition of $U$ and $inj$, the embedding of $G_t$ extends the embedding of $G_{t-1}$, thus we can color the new vertex $v_t$ to $c'(v_t)=c(inj_{G_t}(v_t))$ and it remains true that every vertex $v_i (i\le t)$ of $G_t$ is colored to $c'(v_i)=c(inj_{G_t}(v_i))$.

As the coloring $c$ on $U$ is an offline nonrepetitive coloring, $c$ is a nonrepetitive coloring also on any subgraph of $U$. In particular, during the coloring process at any time it is a nonrepetitive coloring of $inj_{G_t}(G_t)$. As $inj_{G_t}$ injects $G_t$ into $U$ as a subgraph of $U$, the coloring $c'$ of $G_t$ is also a nonrepetitive coloring of $G_t$.

Now we only need to construct $U$ as a sequence of graphs $U_1\subset U_2\subset \dots =U$. This is done by choosing $U_1=K_{k+1}$, the complete graph on $k+1$ vertices, and then $U_i$ is defined recursively by adding a new vertex $v$ for each $k$-clique of $U_{i-1}$ and connecting $v$ to the vertices of the clique. In words, in each step we extend $U_{i-1}$ in all possible ways that the online construction of $k$-trees allows. It is easy to see that there is a natural embedding of all $G_t$ to some $U_s$, $s\le t$, which gives an injection for all $G_t$ into $U$, so that it is naturally incremental in the needed way. As we start with $U_1$ being a $k$-tree and extending $U_{i-1}$ to $U_{i}$ is just a sequence of operations where we add a new vertex $v$ for a $k$-clique and connect it to the vertices of the clique, all the $U_i$'s are $k$-trees by induction, just like we needed.
\end{proof}

\begin{figure}
    \centering
    \subfigure[]{\label{fig:hinfdef}
        \includegraphics[scale=1]{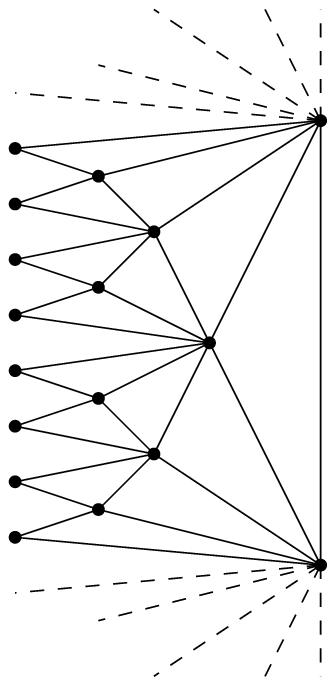}
        }
        \hspace{10mm}
    \subfigure[]{\label{fig:hinfexample}
    		\includegraphics[scale=1]{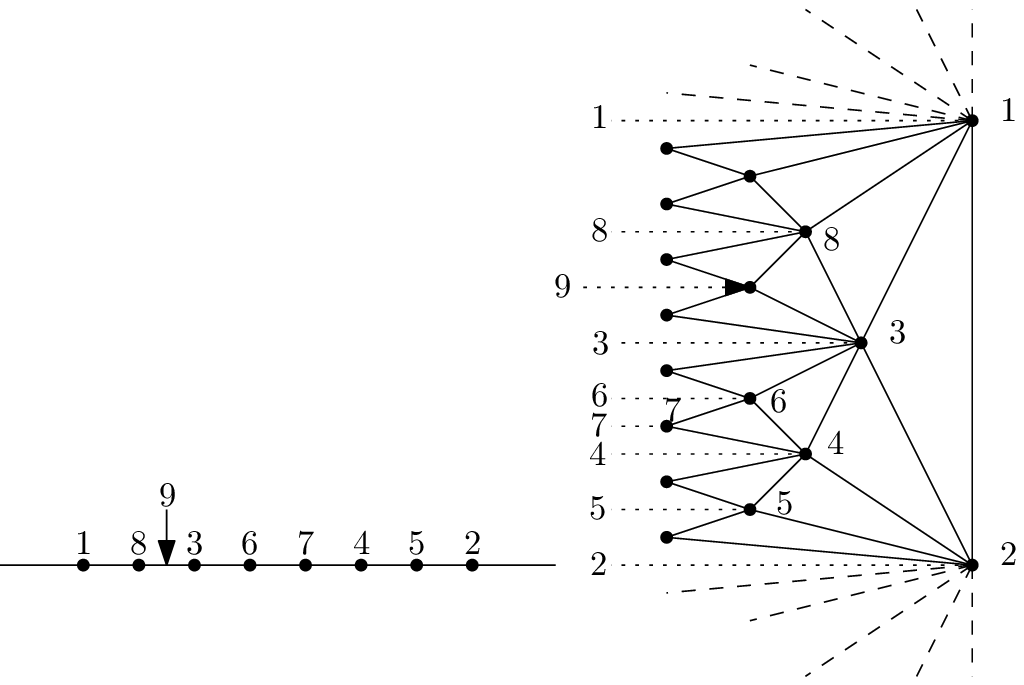}
      }
   \caption{(a) A drawing of a part of $O_5$ and (b) injection of an online path on $9$ vertices to $O$}   
\end{figure}

\begin{proof}[Proof of \tref{onlinepath}]
We build a universal graph $O$ similarly as in the previous proof. 
$O_1$ is a two-vertex graph with an edge.
For $i \ge 0$, the graph $O_{i+1}$ is obtained from $O_i$ by adding some vertices and edges.
So $O_i$ is a subgraph of $O_{i+1}$. The edges of $O_i$ are called {\em old edges} of $O_{i+1}$
and edges in $E(O_{i+1})-E(O_i)$ are called {\em new edges}. The set of new edges of $O_{i+1}$ forms a path.
The two end vertices of the path formed by new edges are called the end vertices of $O_{i+1}$.
Initially, the only edge of $O_1$ is a new edge. 
Assume $O_i$ is given. To construct $O_{i+1}$,  for each new edge $e$ of $O_{i}$, we add a vertex $v_e$ and connect $v_e$ to the two end-vertices of $e$. Also, we append a vertex to each of the two end  vertices of $O_i$. 
It is easy to see that we maintained that the pending edges form a path. In words, $O_{i+1}$ is obtained from $O_i$ by adding a vertex in all possible ways that is allowed by the restriction of an online path, but we proceed only on new edges of $O_i$.

Observe that $O_i$ is an outerplanar graph (see Figure \ref{fig:hinfdef}) and also it is easy to see that $O$ is indeed a universal graph of online paths. The proof of this is left to the reader, defining the injections is similar to the previous proof. The main difference from the previous proof is that in the online process now edges may be deleted as well, thus $G_{t-1}$ is not always a subgraph of $G_t$, nevertheless $G_t$ is always a subgraph of $O$. For an example of embedding an online path see Figure \ref{fig:hinfexample}.  It was proved in \cite{KP} that for any outerplanar graph there exists a nonrepetitive $12$-coloring, thus every $O_i$ also admits a nonrepetitive $12$-coloring, which in turn again using K\"onig's Lemma implies that $O$ admits a nonrepetitive $12$-coloring which finally implies that paths are online nonrepetitively $12$-colorable.
\end{proof}

In the above proof we do not need the full strength of the result of \cite{KP} about outerplanar graphs, in fact we just need that only some specific (on the figures the vertically monotone paths) need to be nonrepetitively colored, which may be useful for possible improvements on the needed number of colors. To make this statement exact, we generalize nonrepetitive colorings to directed trees and more generally to directed graphs (which seems to be an interesting notion in its own right).

\begin{definition}
A vertical path in a rooted tree $T$ is a simple path whose first vertex is a descendant of the last or vice versa. A coloring of a rooted tree $T$ is vertically-nonrepetitive if there is no repetitive sequence among color sequences of vertical paths in $T$.
\end{definition}

This definition was implicitly present in the literature. In the next section we show a result about such coloring of trees. We generalize this further to general directed graphs.

\begin{definition}
A coloring of a directed graph $G$ is directed-nonrepetitive if there is no repetitive sequence among color sequences of  directed paths in $T$. 
\end{definition}

Note that given a non-directed graph $H$ if we direct all its edges in both directions, then the new directed graph $G$ is directed-nonrepetitively colored if and only if under the same coloring $H$ is nonrepetitively colored. Thus this notion generalizes to directed graphs the original notion of nonrepetitive colorings (of undirected graphs).

\begin{definition}
A coloring of a planar graph $O$ embedded in the plane with no horizontal edges is vertically-nonrepetitive if after directing the edges from top to bottom, this directed graph $G$ is directed-nonrepetitive.
\end{definition}

During the rest of this section by a graph we always mean an outerplanar graph together with its embedding.

\begin{problem}\label{prob:op}
What is the minimal constant $s_p$ such that for any embedded planar graph $G$ there exists a vertically-nonrepetitive coloring of $G$ using $s_p$ colors?
What is the minimal constant $s_o$ such that for any embedded outerplanar graph $G$ there exists a vertically-nonrepetitive coloring of $G$ using $s_o$ colors?
\end{problem}

By the result of \cite{KP} $s_o$ is at most $12$ as there exists a $12$-coloring such that there are no repetitive paths at all (and thus there are no vertical repetitive paths). Furthermore, by the very recent result of \cite{thueplanar} $s_p$ is at most $768$ as there exists a $768$-coloring such that there are no repetitive paths at all.

In the next problem by $O$ we mean the same universal graph for online paths defined in the proof of \tref{onlinepath}, embedded in the way suggested by Figure \ref{fig:hinfdef}.

\begin{problem}\label{prob:opspec}
What is the minimal constant $s_o'$ such that there exists a vertically-nonrepetitive coloring of $O$ using $s_o'$ colors?
\end{problem}

Again, $s_o'$ is at most $12$ by the result about outerplanar graphs. 
As we claimed before, there is a strong connection to the online path coloring problem. Namely, if $s_o'$ is a solution to Problem \ref{prob:opspec} then \tref{onlinepath} holds even in the case if we change $12$ to $s_o'$. To see this, one only has to check that any embedding of $G_t$ is a vertical path in $O$ if it is embedded in the way suggested by Figure \ref{fig:hinfexample}.

\section{Online list-colorings} \label{sec:onlinelist}

The definition of online coloring paths, $k$-trees and other graph classes generalizes in a natural way to list-colorings. The list-coloring version seems to be much harder, as already the innocent-looking problem of coloring online nonrepetitively a left-to-right path is non-trivial. Yet, this case is implied by previous research, as we shall see now.

Indeed, from \cite{KM} it follows that lists of size $4$ are enough for such a left-to-right online list-coloring. Their result which we need is about vertically-nonrepetitive list-colorings of trees (in fact they prove a more general lemma):

\begin{lemma} [Kozik and Micek \cite{KM}]\label{lem:kozik}
For any rooted tree $T$ and any list assignment with lists of size $4$ there exists a vertically-nonrepetitive list-coloring of $T$.
\end{lemma}

We note that without the restriction that we consider only vertical paths, there is no such constant. Indeed, by the result of Fiorenzi, Ochem,  Ossona de Mendez and Zhu \cite{thuechoose} for any $l$ there exists a tree which needs lists of size $l$ to guarantee the existence of a nonrepetitive list-coloring.

\lref{kozik} implies the online result by applying it to the (infinite) rooted tree of all possible list assignments. We omit the details of this otherwise straightforward proof, just note that what is needed is that in any step $t$ depending on the previous lists there is only a finite number of practically different possible lists $L_t$:

\begin{corollary}
If all the lists of vertices have size at least $4$, then a left-to-right path is online nonrepetitively list-colorable.
\end{corollary}

For the other online graph classes we investigated in the previous section, it is unknown if such a constant exists.

\begin{problem}
Is there some constant $c$ such that if all the lists have size at least $c$ then there is an online list-coloring algorithm for paths?

Is there some constant $c(k)$ such that if all the lists have size at least $c(k)$ then there is an online list-coloring algorithm for (partial) $k$-trees?
\end{problem}

\subsubsection*{Acknowledgement} This research was partially done while the first author enjoyed the hospitality of Zhejiang Normal University.

\end{document}